\documentclass[11pt]{article}
\usepackage[margin=1in]{geometry}
\usepackage{amsmath,amssymb,amsthm,microtype}
\usepackage{authblk}

\usepackage{hyperref}
\hypersetup{pdftitle={Every rank-$3$ matroid is Z-Rayleigh},
pdfkeywords={matroid, multivariate Tutte polynomial, Rayleigh inequality, negative correlation}}

\newtheorem{theorem}{Theorem}[section]
\newtheorem{lemma}[theorem]{Lemma}
\newtheorem{conjecture}[theorem]{Conjecture}

\newtheorem{proposition}[theorem]{Proposition}

\newtheorem{problem}[theorem]{Problem}
\theoremstyle{remark}

\newcommand{\cl}{\operatorname{cl}}
\newcommand{\Del}{\Delta}
\newcommand{\coefge}{\succeq}

\begin{document}
\title{Every rank-three matroid is $Z$-Rayleigh}
\author{Qinghong Zhao\thanks{Corresponding author. School of Mathematical Sciences,
Huaqiao University, Quanzhou 362021, China.
Email: \texttt{qzhao@hqu.edu.cn}.}
~and~Xiangqian Zhou\thanks{Department of Mathematics and Statistics, Wright State University, Dayton, OH, 45435, USA.
Email: \texttt{xiangqian.zhou@wright.edu}.}}

\date{}
\maketitle
\vspace{-1.5em}
\begin{abstract}  
We prove that every matroid of rank or corank at most three is
$Z$-Rayleigh: its multivariate Tutte polynomial satisfies the Rayleigh
inequalities for every positive weighting and every parameter
$0<q\le1$. This settles a conjecture of Mani and answers the
rank-$3$ spanning-correlation question of Cocks and the correlated-matroid
question of Semple and Welsh. For simple rank-$3$ matroids, we obtain
an explicit decomposition of each spanning-set Rayleigh difference into
the corresponding basis Rayleigh difference, positively weighted
squares, and a polynomial with nonnegative coefficients. A second
identity compares the normalized multivariate Tutte Rayleigh difference
with $(1-q)^2$ times the spanning-set difference, with a remainder of
the form $q(1-q)R(q,x)$, where $R$ has nonnegative coefficients.
The proof uses the bipartite incidence graph of the lines through two
fixed elements and Wagner's theorem that rank-$3$ matroids are
basis-Rayleigh. The rank bound is sharp, since a rank-four matroid can
fail even the basis Rayleigh inequalities.

\medskip
\noindent\textbf{Keywords:} Matroid; multivariate Tutte polynomial;
Rayleigh inequality; negative correlation

\smallskip
\noindent\textbf{2020 Mathematics Subject Classification:} 05B35, 05A20, 60C05.
\end{abstract}

\section{Introduction}\label{sec:intro}
All matroids and graphs considered here are finite. We refer to
Oxley~\cite{O} for matroid terminology. Let $M$ be a matroid on a
ground set $E$, with rank function $r_M$, and put
$r=r(M)=r_M(E)$ and $n=|E|$. The \emph{corank} of $M$ is
$r(M^*)=n-r$, where $M^*$ is its dual. A subset $A\subseteq E$
is independent if $r_M(A)=|A|$, is spanning if $r_M(A)=r$,
and is a basis if both conditions hold. A matroid is \emph{simple}
if it has no loops or parallel elements.

For variables $x=(x_e:e\in E)$ and $A\subseteq E$, write
$x^A=\prod_{a\in A}x_a$, with an empty product equal to one.
A polynomial is \emph{multiaffine} if its degree in each variable
is at most one. For such a polynomial $F$ and distinct $e,f\in E$,
define its \emph{Rayleigh difference} by
\[
 \Del_{e,f}F=(\partial_e F)(\partial_f F)-F(\partial_{e,f}^{2}F),
\]
where
\[
 \partial_e F=\frac{\partial F}{\partial x_e},
 \qquad
 \partial_{e,f}^{2}F=\frac{\partial^{2}F}{\partial x_e\,\partial x_f}.
\]
A nonzero multiaffine polynomial with nonnegative coefficients is
\emph{Rayleigh} if $(\Del_{e,f}F)(x)\ge0$ for all distinct $e,f$ and all
$x>0$, where $x>0$ means $x_a>0$ for every $a\in E$.

This is a weighted negative-correlation condition. Write
$F(x)=\sum_{A\subseteq E}c_Ax^A$ and assign each subset
$A\subseteq E$ probability $\frac{c_Ax^A}{F(x)}$. Since $F$ is
multiaffine,
\[
\Pr(e\in A)=\frac{x_e(\partial_eF)(x)}{F(x)},
\qquad
\Pr(f\in A)=\frac{x_f(\partial_fF)(x)}{F(x)},
\]
and
\[
\Pr(e,f\in A)
=\frac{x_ex_f(\partial_{e,f}^{2}F)(x)}{F(x)}.
\]
Therefore,
\[
\Pr(e\in A)\Pr(f\in A)-\Pr(e,f\in A)
=\frac{x_ex_f}{F(x)^2}(\Del_{e,f}F)(x).
\]
Since $x_ex_f/F(x)^2>0$, the Rayleigh property is equivalent
to pairwise negative correlation under every positive choice
of the weights; see Wagner~\cite[Section~3]{W08}.

The basis, independent-set and spanning-set generating polynomials are
\[
 B_M(x)=\sum_{B\text{ a basis of }M}x^B,\qquad
 I_M(x)=\sum_{A\text{ independent in }M}x^A,\qquad
 S_M(x)=\sum_{r_M(A)=r}x^A,
\]
where all indexing sets are subsets of $E$. Following
Wagner~\cite{W08}, a matroid $M$ is called $B$-Rayleigh (resp. $I$-Rayleigh, $S$-Rayleigh) if the polynomial $B_M(x)$ (resp. $I_M(x)$, $S_M(x)$) is Rayleigh. Semple and Welsh~\cite{SW} call the $I$-Rayleigh property \emph{independence-correlated} and the $S$-Rayleigh property \emph{spanning-correlated}; a matroid is \emph{correlated} if it is both independence-correlated and spanning-correlated. 

The $B$-Rayleigh property extends the classical Rayleigh monotonicity property of electrical networks from spanning trees of graphs to bases of arbitrary matroids. Choe and Wagner~\cite{CW} called a matroid \emph{Rayleigh} if it is $B$-Rayleigh, that is, if $(\Del_{e,f}B_M)(x)\ge0$ for all distinct $e,f\in E(M)$ and all $x>0$. Wagner~\cite{W} proved that every rank-$3$ matroid is Rayleigh. Semple and Welsh~\cite{SW} studied the analogous questions for independent and spanning sets and asked the following.

\begin{problem}[\cite{SW}]
\label{problem}
Is every rank-$3$ matroid correlated?
\end{problem}

Cocks~\cite{C} proved that every rank-$3$ matroid
is $I$-Rayleigh and asked whether every rank-$3$ matroid is
spanning-correlated~\cite[Section~6, Question~(2)]{C}.
A stronger condition concerns the multivariate Tutte polynomial.
Following the normalization used by Sokal~\cite{Sokal} and
Mani~\cite{M}, put
\[
 Z_M(q,x)=\sum_{A\subseteq E}q^{-r_M(A)}x^A.
\]
For fixed $q>0$, this is a multiaffine polynomial in $x$ with positive
coefficients. The matroid $M$ is \emph{$Z$-Rayleigh} if $Z_M(q,x)$
is Rayleigh in $x$ for every $0<q\le1$. All derivatives in a
Rayleigh difference of $Z_M$ are taken with respect to $x$, with
$q$ fixed.  The $Z$-Rayleigh property is invariant under duality; see
Wagner~\cite{W08} and Mani~\cite{M}. 
It is easy to check that $Z_M(1,x)=\prod_{e\in E}(1+x_e)$ and  $(\Del_{e,f}Z_M)(1,x)=0$ for every $\{e,f\}\subseteq E$. 

Both the spanning-set polynomial and the independent-set polynomial may be obtained as limits: 
\begin{equation}\label{spanning-limit}
S_M(x) = \lim _{q \to 0^+}  q^rZ_M(q,x)
 = \lim_{q\to 0^+} \sum_{A\subseteq E}q^{r-r_M(A)}x^A
\end{equation}
and 
\begin{equation}\label{independence-limit}
 I_M(x) = \lim _{q \to 0^+} Z_M(q,qx)= \lim_{q \to 0^+} \sum_{A\subseteq E}q^{|A|-r_M(A)}x^A
\end{equation}
where $qx=(qx_e:e\in E)$. Positive rescaling of a polynomial or its
variables preserves the Rayleigh inequalities, and these inequalities
are preserved by coefficientwise limits of bounded-degree polynomials.
Thus every $Z$-Rayleigh matroid is both $I$-Rayleigh and $S$-Rayleigh;
see also~\cite{W08,M}. The lowest-degree part of $S_M$ is $B_M$, so
the same argument shows that $S$-Rayleigh implies $B$-Rayleigh.

Mani~\cite{M} proved that the rank-$3$ wheel, whirl,
Fano and non-Fano matroids are $Z$-Rayleigh and proposed the following
conjecture.

\begin{conjecture}[\cite{M}]\label{mani}
Every matroid of rank or corank at most three is $Z$-Rayleigh.
\end{conjecture}

In this paper, we prove Conjecture~\ref{mani}. We also obtain a
comparison of the three Rayleigh differences in the simple rank-$3$
case.

\begin{theorem}\label{main}
Every matroid of rank or corank at most three is $Z$-Rayleigh.
Moreover, if $M$ is simple of rank three, then, for all distinct
$e,f\in E(M)$, all $0<q\le1$ and all $x>0$,
\[
 q^6(\Del_{e,f}Z_M)(q,x)
 \ge (1-q)^2(\Del_{e,f}S_M)(x)
 \ge (1-q)^2(\Del_{e,f}B_M)(x)\ge0.
\]
\end{theorem}

In particular, every matroid of rank or corank at most three is
correlated. In rank three, this answers the question of Cocks and
Problem~\ref{problem}, completing the picture provided by the theorems
of Wagner and Cocks. Semple and Welsh's broader question of whether
there is an $I$-Rayleigh matroid that is not $S$-Rayleigh is recorded
as Problem~15.8.13 in Oxley's book~\cite{O}. Our theorem does not
resolve that question in arbitrary rank, but shows that any such
matroid must have both rank and corank at least four.

The rank threshold is sharp. Choe and Wagner~\cite{CW}
discuss the rank-four matroid $S_8$, with its ground set labelled so that
$(\Del_{1,8}B_{S_8})(\mathbf1)=-16$, where $\mathbf1$ is the all-one vector.
For a rank-$r$ matroid,
\[
 (\Del_{e,f}S_M)(tx)
 =t^{2r-2}(\Del_{e,f}B_M)(x)+O(t^{2r-1})
 \qquad(t\longrightarrow0^+).
\]
Consequently, $S_8$ is not $S$-Rayleigh, and
\eqref{spanning-limit} implies that it is not $Z$-Rayleigh.

Recent work of Fang and Ma~\cite{FM} studies further polynomial
positivity conditions that imply Rayleigh inequalities. Their
G\aa rding-polynomial approach applies to spanning-set polynomials of
several classes, including uniform matroids, series--parallel networks,
and all matroids on at most six elements
\cite[Theorem~12.3]{FM}. They also distinguish the Rayleigh and
G\aa rding properties for the Fano matroid
\cite[Proposition~14.10]{FM}. We refer to~\cite{BH,FM} for the
definitions of these polynomial classes. Our proof instead uses
explicit decompositions adapted to rank-$3$ incidence geometry.

Section~\ref{sec:bipartite} recalls the results and notation used in
the proof and establishes the bipartite polynomial decomposition.
In Section~\ref{sec:proof}, we establish the spanning-set comparison and the multivariate Tutte identity,  and then use them to 
prove Theorem~\ref{main}.

\section{Preliminaries and a polynomial decomposition}
\label{sec:bipartite}
We use the following results. The first supplies the
degree-four inequality in our proof. The second is used to
deduce the answer to Problem~\ref{problem}.

\begin{theorem}[\cite{W}]\label{wagner}
For every rank-$3$ matroid $M$, every pair of distinct elements $e,f$
and every $x>0$, $(\Del_{e,f}B_M)(x)\ge0$.
\end{theorem}

\begin{theorem}[\cite{C}]\label{cocks}
For every rank-$3$ matroid $M$, every pair of distinct elements $e,f$
and every $x>0$, $(\Del_{e,f}I_M)(x)\ge0$.
\end{theorem}

In a graph $G$, a \emph{matching} is a set of edges with no common endpoints, and a \emph{vertex
cover} is a set of vertices meeting every edge.

\begin{theorem}[\cite{D}]\label{konig}
In a bipartite graph, the maximum size of a matching equals the
minimum size of a vertex cover.
\end{theorem}

For a polynomial $F$, let $F_d$ denote its homogeneous part of
total degree $d$, and write $F\coefge0$ if all its coefficients
are nonnegative. For a set $A$ of elements or edges, define
\[
P_A=\prod_{a\in A}(1+x_a)=\sum_{T\subseteq A}x^T,
\qquad
\sigma_A=\sum_{a\in A}x_a.
\]
Thus $P_A$ is the generating polynomial of all subsets of $A$,
and $\sigma_A$ is its linear part.

We will use the following identity for Rayleigh differences.
Let $F$ be a multiaffine polynomial, and fix distinct $e,f$.
Grouping its terms according to whether they contain $x_e$
and $x_f$, we may write
\begin{equation}\label{multiaffine}
F=A+x_eB+x_fC+x_ex_fD,
\end{equation}
where $A,B,C,D$ are polynomials that do not involve $x_e$ or $x_f$.
Then
\[
\partial_eF=B+x_fD,\qquad
\partial_fF=C+x_eD,\qquad
\partial_{e,f}^{2}F=D.
\]
Substituting into the definition gives
\begin{equation}\label{bcad}
\Del_{e,f}F
=(\partial_eF)(\partial_fF)-F(\partial_{e,f}^{2}F)
=BC-AD.
\end{equation}
In particular, $\Del_{e,f}F$ is independent of $x_e$ and $x_f$.

We next establish a polynomial decomposition for bipartite graphs,
which will be used to compare the spanning-set and basis Rayleigh
differences. Let $G=(L,R;E)$ be a simple bipartite graph with bipartition
$(L,R)$ and edge set $E$. For each vertex $v$, let $\delta_G(v)$
denote the set of edges incident with $v$. Associate a variable $x_e$ with each edge $e\in E$, and let
\[
P_E=\prod_{e\in E}(1+x_e),
\qquad
\sigma_E=\sum_{e\in E}x_e.
\]
Define
\[
E_G=\sum_{v\in L}\bigl(P_{\delta_G(v)}-1\bigr),
\qquad
F_G=\sum_{v\in R}\bigl(P_{\delta_G(v)}-1\bigr).
\]
Thus $E_G$ and $F_G$ are the generating polynomials of nonempty
edge subsets contained in a star centred in $L$ and $R$,
respectively. Also, let
\[
H_G=\sum_{\substack{\{e,f\}\subseteq E\\
                    \{e,f\}\text{ is a matching}}}x_ex_f
\]
be the generating polynomial of matchings of size two in $G$.
Finally, put
\[
D_G=E_GF_G+(P_E-1)(H_G-\sigma_E).
\]

For each five-edge subset consisting of a four-cycle and a pendant
edge $z$, let $\{e,f\}$ and $\{s,t\}$ be the two pairs of opposite
edges of the cycle. Set
\begin{equation}\label{squares}
 Q_G=\frac12\sum x_z(x_ex_f-x_sx_t)^2,
\end{equation}
where each such edge subset is counted once (an empty sum is zero).

\begin{lemma}\label{bipartite}
For every simple bipartite graph $G$, there is a polynomial $C_G\coefge0$
with no terms of degree less than five such that
\[
 D_G=(D_G)_4+Q_G+C_G.
\]
\end{lemma}
\begin{proof}

If $E=\varnothing$, then $P_E=1$ and $E_G=F_G=0$. Thus $D_G=0$, and the result holds trivially. So we assume that $E \neq \emptyset$. 

The linear parts of $P_E-1,E_G,F_G$ are $\sigma_E$. Two distinct edges
either share a left endpoint, share a right endpoint, or form a matching.
These cases are exclusive because $G$ is simple. Thus
\[
 H_G=(P_E)_2-(E_G)_2-(F_G)_2.
\]
Consequently, $D_G$ has no terms of degree less than four: its
quadratic part is $\sigma_E^2-\sigma_E^2=0$, and its cubic part is
$\sigma_E((E_G)_2+(F_G)_2)+\sigma_E H_G-\sigma_E (P_E)_2=0$. 

Every exponent in $D_G$ is at most two. Consider a monomial
$x_e^2x^A$, where $e\in E$ and $A\subseteq E\setminus\{e\}$, and
let $\eta_A(e)$ be the number of edges $f\in A$ such that $\{e,f\}$
is a matching.
The terms $(P_E-1)H_G$ and $-(P_E-1)\sigma_E$ contribute $\eta_A(e)$ and $-1$,
respectively. The term $E_GF_G$ contributes one precisely when every edge of
$A$ meets $e$, that is, when $\eta_A(e)=0$. Hence the coefficient
of $x_e^2x^A$ in $D_G$ is $0$ when $\eta_A(e)=0$, and $\eta_A(e) -1$ otherwise; hence, it is $\max\{0, \eta_A(e) -1\}$. 

A monomial with exactly two squared variables has coefficient one
if those two edges form a matching, and zero otherwise. Indeed, its only
possible contribution is from selecting those two edges in $H_G$.
Monomials with at least three squared variables have coefficient zero. Thus any monomial containing at least one squared variable has
nonnegative coefficient in $D_G$.

It remains to consider square-free monomials \(x^A\) of degree at least five, since all degree-four terms are contained in \((D_G)_4\). Let $m_2(A)$ be the number of matchings of size two contained in $A$,
and let $b(A)$ count the ordered partitions $A=A_L\mathbin{\dot\cup}A_R$ into nonempty
parts, where $A_L$ is contained in a star centred in $L$ and $A_R$ is
contained in a star centred in $R$. A nonempty star has a unique centre
in a specified bipartition class, so this counts each contribution from
$E_GF_G$ once. Hence the coefficient of $x^A$ in $D_G$ is $m_2(A) - |A| +b(A)$. 

If $A$ contains a matching of size three, its three two-edge subsets
are matchings. Each of the other $|A|-3$ edges forms a matching with
at least one of these three edges. Hence $m_2(A)\ge |A|$.

Otherwise, Theorem~\ref{konig} gives a minimum vertex cover of size at
most two. If it has one vertex, then $m_2(A)=0$ and $b(A)=|A|$:
the part on the opposite side must be a singleton. Suppose the cover
has two vertices on opposite sides, say $u \in L$ and $v \in R$. Let $n_u$ (resp. $n_v$) be the number of edges in $A \setminus \{uv\}$ incident with $u$ (resp. $v$). By the minimality of $\{u,v\}$, $n_u, n_v \geq 1$. Note that
$m_2(A) = n_u n_v$. Moreover, if $uv \in A$, then $|A|=n_u+n_v+1$ and $b(A) =2$; otherwise $|A| = n_u +n_v$ and $b(A) =1$. In both cases, $m_2(A) - |A| +b(A ) =(n_u-1)(n_v-1) \geq 0$.

Now suppose the two covering vertices $u$ and $v$ are on the same side.  Again, let $n_u$ (resp. $n_v$) be the number of edges of $A$ incident with $u$ (resp. $v$). Let $t$ be the number of common neighbors of $u$ and $v$ in the subgraph with edge set $A$. Then $m_2(A) = n_u n_v -t$  and $|A| = n_u +n_v $. 

If $\min\{n_u,n_v\}=1$, say $n_u=1$, then $t=1$ or $t=0$. In the former case, $m_2(A)-|A|+b(A)=n_v-1-(1+n_v)+2=0$, while in the latter case, $m_2(A)-|A|+b(A)=n_v-(1+n_v)+1=0$.

Now assume that $n_u,n_v\ge2$. Then $b(A)=0$, and hence
\[
m_2(A)-|A|+b(A)
=n_un_v-n_u-n_v-t.
\]
Since $t\le\min\{n_u,n_v\}$, this quantity is nonnegative whenever
$\min\{n_u,n_v\}\ge3$. We may therefore assume, without loss of
generality, that $n_u=2\le n_v$. If $n_v\ge4$, then
\[
n_un_v-n_u-n_v-t
=n_v-2-t
\ge n_v-4
\ge0.
\]
Thus the only remaining case is $n_v\le3$. Since
$|A|=n_u+n_v\ge5$, we must have $n_v=3$ and the coefficient is $6 - t - 5 = 1-t$, which is negative only when $t=2$. 
This is precisely a $4$-cycle with one pendant edge, and its
coefficient is $-1$. Thus the only negative coefficients of degree at least five occur in degree five, and their supports are exactly those appearing in \(Q_G\). For each such support, the corresponding summand of \(Q_G\) has negative term $-x_zx_ex_fx_sx_t$,  which has the required coefficient. Its two positive terms are $\frac{1}{2} x_e^2x_f^2x_z$ and  $\frac{1}{2} x_s^2x_t^2x_z$. Each of these monomials has coefficient one in \(D_G\). Moreover, no such positive monomial occurs in two different summands of $Q_G$: a matching of size two determines the unique pair of cross edges that could complete a four-cycle, and the remaining variable specifies the pendant edge. Hence $(D_G)_5-Q_G\coefge0$.
All coefficients of degree at least six are nonnegative. Therefore, taking $C_G=D_G-(D_G)_4-Q_G$ proves the lemma.
\end{proof}

\section{Proof of the main theorem}\label{sec:proof}
First, we establish a polynomial identity for a simple rank-$3$ matroid.
A \emph{line} is a flat of rank two, and $\cl_M$ denotes closure in $M$.

\begin{proposition}\label{certificate}
Let $M$ be a simple rank-$3$ matroid, and let $e,f\in E(M)$ be
distinct. There is a simple bipartite graph $G$ whose edges are indexed
by $E(M)\setminus\cl_M\{e,f\}$ and a polynomial $R_{M,e,f}\coefge0$
with no terms of degree less than five such that
\[
 \Del_{e,f}S_M=\Del_{e,f}B_M+Q_G+R_{M,e,f}.
\]
Here $Q_G$ is defined by \eqref{squares}.
\end{proposition}
\begin{proof}
Let $\ell=\cl_M\{e,f\}$, $T=\ell\setminus\{e,f\}$ and
$U=E(M)\setminus\ell$. Since $M$ has rank three, $U\ne\varnothing$.
Distinct lines of a simple rank-$3$ matroid meet in at most one
point. Therefore, the intersections with $U$ of the lines through $e$ other than $\ell$ give a partition of $U$, denoted by $ \mathcal{U} _e = \{U_i ^ e: 1 \leq i \leq K_e\}$; similarly, lines through $f$ other than $\ell$ also give a partition  of $U$, denoted by $\mathcal{U}_f = \{U_j ^f: 1 \leq j \leq K_f\}$. Note that $|U_i ^e \cap U_j ^f| \leq 1$.  Define a bipartite graph $G$ with edge set $U$, whose two vertex classes are disjoint copies of $\mathcal U_e$ and $\mathcal U_f$. Each $u\in U$ joins the two blocks containing $u$.  The graph $G$ is a simple bipartite graph.  For the graph $G$, we define the polynomials $P_U$, $\sigma_U$, $E_G$, $F_G$, $H_G$, and $D_G$ as in Section~\ref{sec:bipartite}.

For a line $L$, write $L_U=L\cap U$, and define
\[
 H_0=\sum_{L \cap \{e,f\} = \emptyset}(P_{L_U}-1-\sigma_{L_U}),
 \qquad
 G_t=\sum_{\substack{L\ni t\\L\ne\ell}}(P_{L_U}-1)\quad(t\in T).
\]
The sums are over lines of $M$. Terms of $H_0$ with $|L_U|\le1$ are equal to $0$. For two distinct elements $u,v\in U$, let $L$ be the unique line containing them. If $L$ contains $e$ or $f$, then $u$ and $v$ lie in the same block of $\mathcal U_e$ or $\mathcal U_f$, respectively, so the edges $u$ and $v$ of $G$ are adjacent. Conversely, if $L\cap \{e,f\}=\varnothing$, then $u$ and $v$ lie in different blocks of both partitions, and hence form a matching in $G$. Therefore, the quadratic part of $H_0$ is exactly $H_G$. Thus $H_0-H_G\coefge0$
and has no terms of degree less than three. For each $t\in T$, the
sets $L_U$ in the definition of $G_t$ partition $U$. It follows that
$G_t-\sigma_U \coefge0$ and has no terms of degree less than two.

Write $S_M=A+x_eB+x_fC+x_ex_fD$ as in \eqref{multiaffine}. We now find each of $A$, $B$, $C$, and $D$. 

A set containing both $e$ and $f$ spans precisely when it meets $U$,
so $D=P_T(P_U-1)$. A set containing $e$ but not $f$ does not span either
because its other elements lie in $T$, or because they form a nonempty
subset of a line through $e$ other than $\ell$. In this setting, $E_G$ is the generating polynomial of the nonempty subsets of $U$ contained in a line through $e$ other than $\ell$. Therefore, $B=D-E_G$.
By symmetry, $C=D-F_G$. The nonspanning subsets of $T\cup U$ have generating polynomial
\[
 P_T+E_G+F_G-\sigma_U+H_0+\sum_{t\in T}x_tG_t.
\]
The first term counts subsets of $T$. The next three terms count
nonempty subsets of $U$ on a line through $e$ or $f$ other than $\ell$; $-\sigma_U$ corrects for the double count of singletons. 
The term \(H_0\) counts the subsets of \(U\) of size at least two that lie on a line containing neither \(e\) nor \(f\).  A set containing two points of $T$ and
a point of $U$ spans $M$. Thus a nonspanning set meeting both $T$ and
$U$ contains exactly one point of $T$, and the last sum counts all
such sets. Since $A$ is the generating polynomial of all spanning sets contained in $T \cup U$,
\[
\begin{aligned}
A&=P_TP_U-\left(P_T+E_G+F_G-\sigma_U+H_0+
                  \sum_{t\in T}x_tG_t\right)\\
 &=D-E_G-F_G-H_0+\sigma_U-\sum_{t\in T}x_tG_t.
\end{aligned}
\] 

By \eqref{bcad},
\begin{equation}\label{reduction}
 \Del_{e,f}S_M=BC-AD = E_GF_G+P_T(P_U-1)(H_0-\sigma_U + \sum_{t \in T} x_t G_t).
\end{equation}

Recall that $D_G = E_G F_G + (P_U -1 ) (H_G - \sigma _U)$ since $U$ is the edge set of $G$.   By rearranging terms in \eqref{reduction}, we get that  
\[
\Del_{e,f}S_M-D_G=(P_U-1)\bigg[ P_T(H_0-H_G)+(P_T-1)H_G -P_T\sigma_U+\sigma_U+P_T\sum_{t\in T}x_tG_t\bigg].
\]  

Inside the brackets, we add and subtract $P_T\sigma_U\sigma_T$.
After simplifying, we obtain
\begin{equation}\label{nonnegative}
\begin{split}
\Del_{e,f}S_M-D_G=(P_U-1)\bigg[&P_T(H_0-H_G)+(P_T-1)H_G\\
 &+\sigma_U(P_T\sigma_T-P_T+1)\\
 &+P_T\sum_{t\in T}x_t(G_t-\sigma_U)\bigg].
\end{split}
\end{equation} 

Note that
\[
P_T\sigma_T-P_T+1
=\sum_{\varnothing\ne A\subseteq T}(|A|-1)x^A
 +\sum_{t\in T}x_t^2P_{T\setminus\{t\}}.
\]
It follows that $P_T \sigma_T -P_T+1\coefge0$ and has no terms of degree less than two. Each term inside the brackets in
\eqref{nonnegative} has nonnegative coefficients and no terms of
degree less than three. Therefore $\Del_{e,f}S_M-D_G\coefge0$ and has no terms of
degree less than four.

The degree-three part of $S_M$ is $B_M$, and $S_M$ has no terms of
lower degree. Consequently
\[
 \Del_{e,f}B_M=(\Del_{e,f}S_M)_4=(D_G)_4+(\Del_{e,f}S_M-D_G)_4.
\]
Lemma~\ref{bipartite} now gives
\[
 \Del_{e,f}S_M=\Del_{e,f}B_M+Q_G+C_G+
 (\Del_{e,f}S_M-D_G)-(\Del_{e,f}S_M-D_G)_4.
 \]
Since $(\Del_{e,f}S_M-D_G)-(\Del_{e,f}S_M-D_G)_4\coefge0$ and has no terms of degree less than five,
we may take $R_{M,e,f}=C_G+(\Del_{e,f}S_M-D_G)-(\Del_{e,f}S_M-D_G)_4$.
\end{proof}

By Proposition~\ref{certificate} and Theorem~\ref{wagner},  we have that
\begin{equation}\label{simple-spanning}
 (\Del_{e,f}S_M)(x)\ge(\Del_{e,f}B_M)(x)\ge0
\end{equation}
for every simple rank-$3$ matroid $M$, all distinct $e,f$ and all $x>0$.
We now use this inequality to prove the $Z$-Rayleigh property.

For a rank-$3$ matroid $M$, define
\[
 W_M(q,x)=\sum_{A\subseteq E(M)}q^{3-r_M(A)}x^A.
\]
Thus $W_M(q,x)=q^3Z_M(q,x)$ for $q>0$, and $W_M(0,x)=S_M(x)$.
All derivatives below are taken with respect to the variables $x$,
with $q$ fixed.

\begin{proposition}\label{z-certificate}
Let $M$ be a simple rank-$3$ matroid, and let $e,f\in E(M)$ be distinct.
There is a polynomial $\mathcal R_{M,e,f}(q,x)$ with nonnegative
coefficients in $q$ and $x$ such that
\begin{equation}\label{z-comparison}
 (\Del_{e,f}W_M)(q,x)
 =(1-q)^2(\Del_{e,f}S_M)(x)
   +q(1-q)\mathcal R_{M,e,f}(q,x).
\end{equation}
\end{proposition}
\begin{proof}
Use the graph $G$, the sets $T,U$, and the polynomials
$P_T,P_U,\sigma_T,\sigma_U,E_G,F_G,H_0,G_t$ from the proof of
Proposition~\ref{certificate}. Let
\[
 \beta_G=P_U-1+\sigma_U-E_G-F_G.
\]
The constant and linear terms of $\beta_G$ are zero. For
$A\subseteq U$ with $|A|\ge2$, the coefficient of $x^A$ in
$\beta_G$ is one if $A$ is not contained in a star, and is zero
otherwise. Note that $A$ cannot be contained in both such stars,
since $G$ is simple. Therefore,
\begin{equation}\label{beta-positive}
 \beta_G\coefge0.
\end{equation}

We next show that
\begin{equation}\label{H0-beta}
 P_UH_0-\beta_G\coefge0.
\end{equation}
Suppose that $x^A$ has coefficient one in $\beta_G$. Then $A$
contains a matching of size two. Otherwise, any two edges of $A$
meet. Choose two of them, $vw_1$ and $vw_2$. Any edge meeting both
but not incident with $v$ would have to be $w_1w_2$, contrary to
bipartiteness. Thus every edge of $A$ is incident with $v$, which
would place $A$ in a star, a contradiction. Let $\{g,h\}\subseteq A$
be a matching. Since the quadratic part of $H_0$ is $H_G$, the
monomial $x_gx_h$ has coefficient one in $H_0$. Also,
$x^{A\setminus\{g,h\}}$ occurs in $P_U$. Hence the coefficient of
$x^A$ in $P_UH_0$ is at least one. All coefficients of $P_UH_0$
are nonnegative, so \eqref{H0-beta} follows.

Recall that we established the following inequalities in the proof of Proposition~\ref{certificate}.
\[
 G_t-\sigma_U\coefge0\quad(t\in T),
 \qquad P_T\sigma_T-P_T+1\coefge0.
\]

Write $W_M=A_q+x_eB_q+x_fC_q+x_ex_fD_q$, where $A_q$, $B_q$, $C_q$, and $D_q$ do not involve $x_e$ or $x_f$. We now
find these four polynomials. A set containing both $e$ and $f$
has rank three precisely when it meets $U$; otherwise it has rank
two. Therefore, $D_q=P_T(P_U-1+q)$.

We next determine $B_q$ by comparing it with $D_q$. For each
$S\subseteq T\cup U$, the monomial $x^S$ in $D_q$ corresponds to
$S\cup\{e,f\}$, while the same monomial in $B_q$ corresponds to
$S\cup\{e\}$. Thus we need to determine when removing $f$ changes
the rank, and hence the coefficient. Suppose that $S$ is a nonempty subset of $U$ contained in a line
through $e$ other than $\ell$. Then $S\cup\{e,f\}$ has rank three,
while $S\cup\{e\}$ has rank two. Hence the coefficient of $x^S$
changes from $1$ to $q$. Since these subsets are counted exactly
by $E_G$, their total contribution to $B_q-D_q$ is
$(q-1)E_G$.

If $S=\varnothing$, the two sets are $\{e,f\}$ and $\{e\}$,
which have ranks two and one, respectively. Thus the constant
term changes from $q$ to $q^2$, giving the additional correction
$q^2-q$.

All other choices have unchanged rank. Indeed, if $S$ is a
nonempty subset of $T$, both sets have rank two. If $S$ meets
both $T$ and $U$, both sets span $M$. Finally, if $S\subseteq U$
is nonempty and is not contained in a line through $e$, then
$S\cup\{e\}$ already spans $M$.

Therefore,
\[
\begin{aligned}
B_q
&=D_q+(q-1)E_G+(q^2-q)\\
&=D_q-(1-q)E_G-q(1-q)\\
&=D_q-(1-q)(E_G+q).
\end{aligned}
\]

By symmetry, $C_q= D_q -(1-q) (F_G +q)$. 

Recall that $A$ is the generating polynomial of the spanning
subsets of $T\cup U$. Since $M$ is simple, the empty set is its
only rank-zero subset, and its rank-one subsets are the singletons.
Thus the rank-two subsets of $T\cup U$ have generating polynomial $P_TP_U-A-1-\sigma_U-\sigma_T$. 
It follows that
\[
 A_q=A+q(P_TP_U-A-1-\sigma_U-\sigma_T) +q^2(\sigma_U+\sigma_T)+q^3.
\]
From the proof of Proposition~\ref{certificate}, we know that 
\[
A=P_TP_U-\left(P_T+E_G+F_G-\sigma_U+H_0+\sum_{t\in T}x_tG_t\right). 
\] 
By substituting $A$ in $A_q$ and rearranging the terms, we get that 
\begin{equation}\label{Aq}
 A_q=D_q-(1-q)\bigg[E_G+F_G+H_0-\sigma_U +\sum_{t\in T}x_tG_t +q(\sigma_U+\sigma_T+1+q)\bigg].
\end{equation}

By \eqref{bcad},  $\Del_{e,f}W_M=B_q C_q - A_q D_q$; then we plug in each of $A_q,\ B_q ,\ C_q$ and $\ D_q$ and rearrange: 

 \begin{equation}\label{Delta-W-first}
 \Del_{e,f}W_M=(1-q)\bigg[  (1-q)(E_G+q)(F_G+q)+P_T(P_U-1+q)\bigg(H_0-\sigma_U +\sum_{t\in T}x_tG_t +q(\sigma_U+\sigma_T-1+q)\bigg)\bigg].
\end{equation}

Recall Equation \eqref{reduction} gives  $\Del_{e,f}S_M= E_GF_G+P_T(P_U-1)(H_0-\sigma_U + \sum_{t \in T} x_t G_t)$. Subtracting $(1-q)^2$ times \eqref{reduction} from \eqref{Delta-W-first}, we obtain
\begin{equation}\label{raw-remainder}
\begin{split}
 &\Del_{e,f}W_M-(1-q)^2(\Del_{e,f}S_M) =q(1-q)\bigg[
 P_TP_U\bigg(H_0-\sigma_U+\sum_{t\in T}x_tG_t\bigg)\\
 &\hspace{29mm}+P_T(P_U-1+q)(\sigma_U+\sigma_T-1+q)+(1-q)(E_G+F_G+q)\bigg].
\end{split}
\end{equation}

We now rearrange the expression inside the brackets. Recall that $ \sum_{t\in T}x_tG_t =\sigma_U\sigma_T+ \sum_{t\in T}x_t(G_t-\sigma_U)$ and $E_G+F_G=P_U-1+\sigma_U-\beta_G$.  By substitution and simplification, the quantity inside the bracket is equal to the following: 
\[
\begin{aligned}
 &P_TP_UH_0-(1-q)\beta_G+P_T\sigma_T(P_U\sigma_U+P_U-1+q)\\
 &+(P_T-1)\bigl[(P_U-1+q)(\sigma_U+q)
               -(P_U\sigma_U+P_U-1+q)\bigr]+P_TP_U\sum_{t\in T}x_t(G_t-\sigma_U).
\end{aligned}
\]

Now we write $P_TP_UH_0=P_UH_0+(P_T-1)P_UH_0$ and rearrange the middle three terms. We get 
\begin{equation}\label{positive-remainder}
\begin{split}
 \mathcal R_{M,e,f}(q,x)={}&(P_UH_0-\beta_G)+q\beta_G\\
 &+(P_T-1)\bigl[P_UH_0+(P_U-1+q)(\sigma_U+q)\bigr]\\
 &+(P_T\sigma_T-P_T+1)(P_U\sigma_U+P_U-1+q)\\
 &+P_TP_U\sum_{t\in T}x_t(G_t-\sigma_U).
\end{split}
\end{equation}
The first two terms have nonnegative coefficients by
\eqref{beta-positive} and \eqref{H0-beta}. Each factor in the
remaining terms also has nonnegative coefficients, since
$P_T-1\coefge0$, $P_U-1\coefge0$,
$P_T\sigma_T-P_T+1\coefge0$, and $G_t-\sigma_U\coefge0$ for
all $t\in T$. Therefore $\mathcal R_{M,e,f}(q,x)\coefge0$ as a
polynomial in $q$ and $x$. Together with \eqref{raw-remainder},
this proves \eqref{z-comparison}.
\end{proof}

\begin{proof}[Proof of Theorem~\ref{main}]
First, suppose that $M$ is simple of rank three. Since
$W_M=q^3Z_M$, we have
$\Del_{e,f}W_M=q^6(\Del_{e,f}Z_M)$. By
Proposition~\ref{z-certificate} and \eqref{simple-spanning},
\[
\begin{aligned}
 q^6(\Del_{e,f}Z_M)(q,x)
 &\ge(1-q)^2(\Del_{e,f}S_M)(x)\\
 &\ge(1-q)^2(\Del_{e,f}B_M)(x)\ge0
\end{aligned}
\]
for $0<q\le1$ and $x>0$. Thus the theorem holds in this case.

Now let $M$ be any rank-$3$ matroid, and let $\bar M$ be
its simplification. We identify the elements of $\bar M$
with the nonloop parallel classes of $M$. Let $L$ be the set of
loops of $M$, and for each nonloop parallel class $D$, put
\[
 y_D=P_D-1=\prod_{i\in D}(1+x_i)-1.
\]
The rank of a subset of $E(M)$ depends only on the nonloop
parallel classes that it meets. The polynomial $y_D$ counts the
nonempty choices from $D$, while $P_L$ counts all choices of
loops. Therefore,
\begin{equation}\label{simplification-identity}
 Z_M(q,x)=P_LZ_{\bar M}(q,y).
\end{equation}

Let $i,j\in E(M)$ be distinct. If $i$ or $j$ is a loop, then
its factor $1+x_i$ or $1+x_j$ makes the Rayleigh difference zero.
Suppose that $i$ and $j$ belong to distinct nonloop parallel
classes $D$ and $D'$. By the chain rule,
\[
 (\Del_{i,j}Z_M)(q,x)
 =P_L^2(\partial_i y_D)(\partial_j y_{D'})
       (\Del_{D,D'}Z_{\bar M})(q,y)\ge0.
\]
Indeed, all coordinates of $y$ and both derivatives are positive
when $x>0$, and $\bar M$ is a simple rank-$3$ matroid.

It remains to consider $i,j$ in the same nonloop parallel class
$D$. Write
\[
 Z_{\bar M}(q,y)=A+y_DB,
\]
where $A$ and $B$ do not involve $y_D$. For
$X\subseteq E(\bar M)\setminus\{D\}$, the coefficient of
$y^X$ in $B-A$ is
\[
 q^{-r_{\bar M}(X\cup\{D\})}
 -q^{-r_{\bar M}(X)}\ge0,
\]
since $r_{\bar M}(X\cup\{D\})\ge r_{\bar M}(X)$ and
$0<q\le1$. Thus $B-A$ has nonnegative coefficients in the remaining
variables $y$ for each fixed $0<q\le1$. Note that
\[
 y_D=P_{D\setminus\{i,j\}}(1+x_i)(1+x_j)-1.
\]
Substituting this into \eqref{simplification-identity} and applying
\eqref{bcad}, we obtain
\[
 (\Del_{i,j}Z_M)(q,x)
 =P_L^2P_{D\setminus\{i,j\}}B(B-A)\ge0.
\]
Here $B$ also has nonnegative coefficients. Therefore every
rank-$3$ matroid is $Z$-Rayleigh.

If $r(M)<3$, add $3-r(M)$ coloops to obtain a rank-$3$ matroid
$N$. Then
\[
 Z_N(q,x)=Z_M(q,x)\prod_{c\text{ a new coloop}}(1+x_c/q).
\]
For any two elements of $M$, the corresponding Rayleigh difference
is multiplied by the square of the displayed positive product.
Since $N$ is $Z$-Rayleigh, so is $M$.

Since the $Z$-Rayleigh property is dual invariant, every matroid with corank at most $3$ is also $Z$-Rayleigh. 
\end{proof}

By \eqref{spanning-limit}, Theorem~\ref{main} implies that every
rank-$3$ matroid is $S$-Rayleigh. By Theorem~\ref{cocks}, it is
also $I$-Rayleigh. Therefore, we establish that every rank-3 matroid is correlated, providing affirmative answers to both Problem~\ref{problem}  and the spanning-correlation question of Cocks. 

\section*{AI Declaration}
The authors used ChatGPT Pro to assist with checking mathematical arguments, revising the text and locating references. The authors take full responsibility for all results, proofs, citations and the final content of this manuscript.

\end{document}